\documentclass[11pt]{amsart}
\usepackage{amssymb}
\usepackage{fullpage}
\usepackage{url}
\usepackage{enumitem}

\newcommand{\corref}[1]{Corollary~\ref{#1}}
\newcommand{\thmref}[1]{Theorem~\ref{#1}}
\newcommand{\lemref}[1]{Lemma~\ref{#1}}
\newcommand{\propref}[1]{Proposition~\ref{#1}}

\newcommand{\tabref}[1]{Table~\ref{#1}}
\newtheorem{thm}{Theorem}[section] %to number theorems, etc. within sections
\newtheorem{prop}[thm]{Proposition}
\newtheorem{cor}[thm]{Corollary}
\newtheorem{lem}[thm]{Lemma}

\theoremstyle{remark}

\newtheorem*{ack}{Acknowledgment}

\numberwithin{equation}{section}  %to number equations within sections

\DeclareMathOperator{\ch}{\operatorname{char}}

\newcommand{\Bol}{\textsc{Bol}} % Bol identity
\newcommand{\Mfg}{\textsc{Mfg}} % Moufang identity
\newcommand{\RIP}{\textsc{Rip}} % right inverse property
\newcommand{\RAP}{\textsc{Rap}} % right alternative property
\newcommand{\SAIP}{\textsc{Saip}} % semiautomorphic inverse property

\begin{document}

\title{When is a Bol loop Moufang?}

\author{Orin Chein}
\address{Temple University \\ Philadelphia, PA 19122 U.S.A.}
\email{\url{orin@math.temple.edu}}
\author{Edgar G. Goodaire}
\address{Memorial University of Newfoundland\\
St. John's, Newfoundland\\ Canada A1C 5S7}
\email{\url{edgar@mun.ca}}
\author{Michael Kinyon}
\address{University of Denver\\Denver, CO 80208 U.S.A.}
\email{\url{mkinyon@math.du.edu}}
\urladdr{\url{www.math.du.edu/~mkinyon}}

\date{\today}

\begin{abstract}
There are a number of identities which, if satisfied
by a Bol loop, imply that the loop is actually Moufang.
In this paper we show that in a number of cases,
the Moufang identity is also forced not by a single
identity, but by giving elements a choice of equations
to satisfy.
\end{abstract}

\maketitle

\section{Background}

A loop is \emph{(right) Bol} if it satisfies
\begin{align}
(xy\cdot z)y = x(yz\cdot y) & \qquad\text{the right Bol identity} \tag{\Bol} \\
\intertext{and \emph{Moufang} if it satisfies any of several equivalent identities including}
(xy\cdot z)y = x(y\cdot zy) & \qquad\text{the right Moufang identity.} \tag{\Mfg}
\end{align}
Despite the similarity in these identities, there are striking
differences between Bol loops and Moufang loops.  For example, while
a Moufang loop is \emph{diassociative} in the sense that any two elements
generate a group, a Bol loop cannot satisfy even rather mild instances
of diassociativity without being Moufang. The following well-known result
lists a few of these.

\begin{prop}[\cite{Robinson:66,Pflugfelder:90}]\label{thm1}
Let $L$ be a (right) Bol loop.  Then any one of the following
identities implies that $L$ is Moufang.
\begin{center}
\begin{tabular}{lll}
(1) & (flexibility)                      & $(xy)x= x(yx)$ \\
(2) & (left alternative law)             & $x(xy)=x^2y$ \\
(3) & (left inverse property)            &  $x^{-1}(xy)=y$ \\
(4) & (commutativity)                    & $xy=yx$ \\
(5) & (antiautomorphic inverse property) & $(xy)^{-1}=y^{-1}x^{-1}$ \\
(6) & (weak inverse property, form 1)     & $x(yx)^{-1}=y^{-1}$ \\
(7) & (weak inverse property, form 2)     & $(xy)^{-1}x=y^{-1}$. \\
\end{tabular}
\end{center}
\end{prop}

The Moufang identity is sometimes forced on a Bol loop not by any
of the two-variable identities that appear in \propref{thm1},
but rather by giving a pair of elements a choice of two equations
to satisfy.   \tabref{tab1} shows the equations we have in
mind, all related to the identities of \propref{thm1}
but interpreted as equations for $x$ and $y$, or $y$ and $x$,
and suggests $\binom{13}{2}$ possible theorems.  In fact, the
number is considerably smaller.

\begin{table}[ht]
\small
\renewcommand{\arraystretch}{1.2}
\begin{tabular}{lll}
commutativity & $xy=yx$ \\
flexibility & $(xy)x= x(yx)$ & $(yx)y = y(xy)$ \\
left alternativity & $x(xy)=x^2y$ & $y(yx)=y^2x$ \\
left inverse property & $x^{-1}(xy)=y$ & $y^{-1}(yx)= x$ \\
antiautomorphic inverse property & $(xy)^{-1}=y^{-1}x^{-1}$ & $(yx)^{-1}=x^{-1}y^{-1}$\\
weak inverse property, form 1 & $x(yx)^{-1}=y^{-1}$ & $y(xy)^{-1}=x^{-1}$ \\
weak inverse property, form 2 & $(xy)^{-1}x=y^{-1}$ & $(yx)^{-1}y = x^{-1}$
\end{tabular}
\par\medskip
\caption{}
\label{tab1}
\end{table}

Recall that Bol loops satisfy the \emph{right inverse property}
$(xy)y^{-1} = x$. This property shows that $x(yx)^{-1}=y^{-1}$ if and
only if $y^{-1}(yx)=x$
and $y(xy)^{-1}=x^{-1}$ if and only if $x^{-1}(xy)=y$ so we can remove
the first weak inverse property equations, because they are equivalent to those
corresponding to the left inverse property.  Similarly, the equations
corresponding to the second weak inverse property are equivalent to those
corresponding to the antiautomorphic inverse property. Thus all the
weak inverse property equations are redundant. Also, as we will show in
Corollary \corref{cor:1}, the left alternative equations are equivalent
to flexibility, and hence they are also redundant.

In addition, there is obvious symmetry in the variables $x$ and $y$. For
instance, the condition ``for all $x,y$, either
$(xy)x = x(yx)$ or $x^{-1}(xy) = y$'' is obviously equivalent to
``for all $x,y$, either $(yx)y = y(xy)$ or $y^{-1}(yx) = x$.''

The upshot of all these reductions is that we are left with twelve pairs of
equations to consider. At least one pair of equations does not force a Bol
loop to be Moufang. The Bol loop of order $8$ denoted $B_8(\Pi_5)$ by Burn
\cite{Burn:78} has the following property: For every $x,y$, either $xy = yx$
or $(xy)^{-1}=y^{-1}x^{-1}$.

There is one case we have not been able to resolve. Specifically, suppose $L$
is a Bol loop and that for each $x,y\in L$,
\[
\text{either}\qquad (xy)^{-1}=y^{-1}x^{-1}
\qquad\text{or}\qquad (yx)^{-1}=x^{-1}y^{-1}\,.
\]
Need $L$ be Moufang?  We do not know, but see Sections \ref{sec:auto} and \ref{sec:srar}.

The precise statement of our findings for the remaining ten cases follows.

\begin{thm}\label{thm:main} Let $L$ be a Bol loop and suppose that
one of the following situations occurs: For any $x,y\in L$,
\begin{center}
\begin{tabular}{rlclc}
(1) & either & $xy=yx$        & or & $x^{-1}(xy)=y$.          \\
(2) & either & $xy=yx$        & or & $(xy)x = x(yx)$.         \\
(3) & either & $x^{-1}(xy)=y$ & or & $(xy)^{-1}=y^{-1}x^{-1}$.\\
(4) & either & $(xy)x=x(yx)$  & or & $(yx)^{-1}=x^{-1}y^{-1}$.\\
(5) & either & $x^{-1}(xy)=y$ & or & $y^{-1}(yx)=x$.          \\
(6) & either & $(xy)x=x(yx)$  & or & $(xy)^{-1}=y^{-1}x^{-1}$.\\
(7) & either & $x^{-1}(xy)=y$ & or & $(yx)^{-1}=x^{-1}y^{-1}$.\\
(8) & either & $(xy)x=x(yx)$  & or & $y^{-1}(yx)=x$.          \\
(9) & either & $(xy)x=x(yx)$ & or & $x^{-1}(xy)=y$.           \\
(10) & either & $(xy)x=x(yx)$  & or & $(yx)y=y(xy)$.
\end{tabular}
\end{center}
Then $L$ is a Moufang loop.
\end{thm}

\begin{ack}
We are pleased to acknowledge the assistance of the automated deduction tool
\textsc{Prover9} \cite{Prover9}.
\end{ack}

\section{Preliminaries}\label{sec:preliminaries}

In this section, we record a few facts about Bol loops we will need later.
In a loop $L$, the \emph{right translation} permutations
$R(a) : L \to L$ are defined, as usual, by $x R(a) = xa$ for all $x\in L$.
An advantage of working with right translations is that they associate,
\emph{i.e.}, $[R(x)R(y)]R(z) = R(x)[R(y)R(z)]$.
The Bol identity itself can be expressed succinctly in terms of translations by
\[
R(x)R(y)R(x) = R((xy)x)\,. \tag{\Bol}
\]
We will reference either the equational or translational form of the Bol
identity by (\Bol).

We will use frequently the fact that a Bol loop satisfies the
\emph{right inverse property}
\[
(yx)x^{-1} = y
\quad\text{or equivalently,}\quad
R(x)^{-1} = R(x^{-1})\,, \tag{\RIP}
\]
and the \emph{right alternative property}
\[
(yx)x = yx^2
\quad\text{or equivalently,}\quad
R(x)^2 = R(x^2)\,.  \tag{\RAP}
\]
Perhaps less familiar, but still very useful is the \emph{semiautomorphic inverse property}
\[
[(xy)x]^{-1} = (x^{-1}y^{-1})x^{-1}
\tag{\SAIP}
\]
see \cite[Thm.~IV.6.12]{Pflugfelder:90}, \cite[Thm.~2.4]{Robinson:66}.

\begin{lem}\label{lem:solution}
If $L$ is a Bol loop, then for $x,y\in L$, the unique solution $z$ to
the equation $yz = x$ is
\begin{equation}
\label{eqn:soln}
z = [x(yx)^{-1}]x = [y^{-1}(xy)]y^{-1}\,.
\end{equation}
\end{lem}

\begin{proof}
Using (\Bol) and (\RIP), it is easy to see that both
sides are solutions. Their equality follows since $L$ is a loop.
\end{proof}

\begin{lem}\label{lip-equiv}
Let $L$ be a Bol loop. For each $x,y\in L$, the following are equivalent:
\begin{enumerate}
\item\qquad $x^{-1}(xy) = y$
\item\qquad $x(x^{-1}y) = y$
\item\qquad $(x^{-1}y)x = x^{-1}(yx)$
\item\qquad $(xy)x^{-1} = x(yx^{-1})$
\end{enumerate}
\end{lem}

\begin{proof}
First note that (1) is equivalent to
\begin{equation}
\label{eqn:lip-tmp1}
x^{-1} = y(xy)^{-1}
\end{equation}
by the right inverse property.

Now if (1) holds, then by \eqref{eqn:lip-tmp1} and (\Bol), we have
$x(x^{-1}y)= x\{[y(xy)^{-1}]y\} = [(xy)(xy)^{-1}]y = y$. The converse
follows from switching the roles of $x$ and $x^{-1}$. Thus (1) is equivalent
to (2).

By \lemref{lem:solution} (with the roles of $x$ and $y$ reversed)
and the right inverse property,
\begin{equation}
\label{eqn:lip-tmp2}
x^{-1}(yx) = \{ [y(xy)^{-1}]y\}x\,.
\end{equation}
If (1) holds, then applying \eqref{eqn:lip-tmp1} to
\eqref{eqn:lip-tmp2}, we get $x^{-1}(yx) = (x^{-1}y)x$, that is, (3) holds. Conversely,
if (3) holds, then \eqref{eqn:lip-tmp2} becomes $(x^{-1}y)x = \{ [y(xy)^{-1}]y\}x$.
Canceling $x$ and then $y$ on the right gives $x^{-1} = y(xy)^{-1}$. Thus (1)
is equivalent to (3).

Finally, the equivalence of (2) and (4) follows from the equivalence of (1) and
(3) by reversing the roles of $x$ and $x^{-1}$.
\end{proof}

The standard \emph{right inner mappings}
$R(x,y) : L \to L$ are defined by $R(x,y) = R(x)R(y)R(xy)^{-1}$.

\begin{lem}
\label{lem.5}
Let $L$ be a Bol loop, let $x_1,\ldots,x_n\in L$, and let
$\varphi : L \to L$ be the permutation defined by
\[
\varphi = R(x_n) R(x_{n-1}) \cdots R(x_1) R(a)^{-1}\,,
\]
where $a = 1 R(x_n) R(x_{n-1}) \cdots R(x_1)$. Then
\[
\varphi^{-1} = R(x_1) R(x_2) \cdots R(x_n) R(b)^{-1}\,,
\]
where $b = 1 R(x_1) R(x_2) \cdots R(x_n)$.
In particular,
\begin{equation}
\label{eqn:rim-inv}
R(x,y)^{-1} = R(y,x)\,.
\end{equation}
\end{lem}

\begin{proof}
First, note that $1 \varphi = 1$. Now
set $a_0 = a^{-1}$ and $a_i = (x_i a_{i-1}) x_i$ for $i > 0$.
Then $\varphi = R(x_n) \cdots R(x_1) R(a_0)$ using (\RIP).
Using successive applications of (\Bol) we have
\begin{align*}
\varphi R(x_1) \cdots R(x_n)
&= R(x_n) R(x_{n-1}) \cdots R(x_1) R(a_0) R(x_1) \cdots R(x_{n-1}) R(x_n) \\
&= R(x_n) R(x_{n-1}) \cdots R( a_1 ) \cdots R(x_{n-1}) R(x_n)
= \cdots
= R( a_n )\,.
\end{align*}
Now
\[
a_n = 1 R(a_n) = 1 \varphi R(x_1)\cdots R(x_n)
= 1 R(x_1)\cdots R(x_n) = b\,,
\]
and so
\[
\varphi R(x_1) \cdots R(x_n) R(b)^{-1} = \mathrm{id}_L\,,
\]
as desired.

The remaining assertion is just the special case $n = 2$.
\end{proof}

The special case \eqref{eqn:rim-inv} can be found in \cite{Kiechle:02}
(for left Bol loops). We will need the more general situation of
\lemref{lem.5} later. We can also cite \cite{Kiechle:02} as a reference
for the next result, but since it is of such evident interest, we feel
it appropriate to include a short proof here.

\begin{lem}\label{lem1}
Let $L$ be a Bol loop. Then for all $x,y\in L$,
\begin{align}
R(x,y) &= R(yx,y)\,, \label{eqn:rim-ls} \\
R(x,y) &= R(x,xy)\,,  \label{eqn:rim-rs} \\
R(x,y) &= R(y^{-1},xy)\,. \label{eqn:rim-s3}
\end{align}
\end{lem}

\begin{proof}
Using (\Bol),
\[
R((yx)y)^{-1} =
[ R(y)R(x)R(y) ]^{-1}
= R(y)^{-1} R(x)^{-1} R(y)^{-1}\,.
\]
Thus
\begin{alignat*}{2}
R(yx,y)
&= R(yx)R(y)R((yx)y)^{-1}
&&= R(yx)R(y) R(y)^{-1} R(x)^{-1} R(y)^{-1} \\
&= R(yx) R(x)^{-1} R(y)^{-1}
&&= [ R(y) R(x) R(yx)^{-1} ]^{-1} \\
&= R(y,x)^{-1}
&&= R(x,y)\,,
\end{alignat*}
using \eqref{eqn:rim-inv} in the last equality.

For \eqref{eqn:rim-rs}, we use \eqref{eqn:rim-inv}, \eqref{eqn:rim-ls}
and \eqref{eqn:rim-inv} again:
$R(x,y) = R(y,x)^{-1} = R(xy,x)^{-1} = R(x,xy)$.

Finally, for \eqref{eqn:rim-s3}, we compute
$R(x,y) = R(x,xy) = R((xy)y^{-1},xy) = R(y^{-1},xy)$ using
\eqref{eqn:rim-rs}, (\RIP) and \eqref{eqn:rim-ls}.
\end{proof}

\begin{cor}\label{cor:1} If $x$ and $y$ are elements of a Bol loop,
then
\[
x(xy)=x^2y
\qquad\text{if and only if}\qquad
(xy)x=x(yx)\,.
\]
\end{cor}
\begin{proof}
We have $x(xy)=(xx)y$ if and only if $xR(x)R(y)R(xy)^{-1}=x$, that
is, if and only if $xR(x,y)=x$.  By \eqref{eqn:rim-inv}, this occurs if and only
if $x=xR(y,x)=xR(y)R(x)R(yx)^{-1}$, that is, if and only if
$(xy)x=x(yx)$.
\end{proof}

\begin{lem}
\label{lem:preaip}
For all $x,y$ in a Bol loop $L$,
\begin{equation}
\label{eqn:preaip}
\{ (xy)[(x^{-1}y^{-1})(xy)]^{-1}\} (xy) = (yx^2)y\,.
\end{equation}
\end{lem}

\begin{proof}
Take the left side of \eqref{eqn:preaip} and multiply it on the left
by $x^{-1}y^{-1}$. Using (\Bol), this reduces to $xy$. On other side,
we have $(x^{-1}y^{-1})[(yx^2)y] = (\{[(x^{-1}y^{-1})y]x\}x)y =
xy$ using (\Bol) and (\RAP), then (\RIP) and then simplifying. Thus both sides
of \eqref{eqn:preaip} are solutions $z$ of the equation $(x^{-1}y^{-1})z = xy$.
Since $L$ is a loop, the two sides are then equal.
\end{proof}

\begin{lem}
\label{lem:decomp}
For all $x,y$ in a Bol loop $L$,
\begin{equation}
\label{eqn:decomp}
R(x^{-1}y^{-1},xy) = R(y^{-1},x^{-1}) R(x,y)\,.
\end{equation}
\end{lem}

\begin{proof}
We compute $R(x^{-1}y^{-1},xy) R(x,y)^{-1} =$
\begin{alignat*}{4}
&= R(x^{-1}y^{-1}) \underbrace{R(xy) R((x^{-1}y^{-1})(xy))^{-1} R(xy)} R(y)^{-1} R(x)^{-1} && \\
&= R(x^{-1}y^{-1}) R(\underbrace{\{(xy)[(x^{-1}y^{-1})(xy)]^{-1}\}(xy)}) R(y)^{-1} R(x)^{-1} &&
\text{by (\RIP), (\Bol)} \\
&= R(x^{-1}y^{-1}) \underbrace{R((yx^2)y)} R(y)^{-1} R(x)^{-1} && \text{by \eqref{eqn:preaip}} \\
&= R(x^{-1}y^{-1}) R(y) R(x)^2 R(y) R(y)^{-1} R(x)^{-1} && \text{by (\Bol), (\RAP)} \\
&= R(x^{-1}y^{-1}) R(y) R(x) &&  \\
&= [ R(x^{-1}) R(y^{-1}) R((x^{-1}y^{-1})^{-1}) ]^{-1} && \text{by (\RIP)} \\
&= R(x^{-1},y^{-1})^{-1} && \\
&= R(y^{-1},x^{-1}) && \text{by \eqref{eqn:rim-inv}}\,.
\end{alignat*}
Thus $R(x^{-1}y^{-1},xy) = R(y^{-1},x^{-1}) R(x,y)$, as claimed.
\end{proof}

\begin{lem}
\label{lem:turnsquare}
For all $x,y$ in a Bol loop $L$,
\begin{equation}
\label{eqn:turnsquare}
R((xy)^2) = R(y)R(x)R([(y^{-1}x^{-1})(yx)]^{-1})R(x)R(y)\,.
\end{equation}
\end{lem}

\begin{proof}
We compute
\begin{alignat*}{4}
(xy)^2 &= (xy)R(xy) R(y)^{-1} R(x)^{-1} R(x) R(y) && \\
&= (xy)R(x,y)^{-1} R(x) R(y) &&  \\
&= (xy)R(y,x) R(x) R(y) && \text{by \eqref{eqn:rim-inv}} \\
&= [(xy^2)x] R(yx)^{-1} R(x) R(y) && \text{by (\RAP)} \\
&= [\{(yx)[(y^{-1}x^{-1})(yx)]^{-1}\}(yx)] R(yx)^{-1} R(x) R(y) && \text{by \eqref{eqn:preaip}} \\
&= \{(yx)[(y^{-1}x^{-1})(yx)]^{-1}\} R(x) R(y) &&  \\
&= [y(\{x[(y^{-1}x^{-1})(yx)]^{-1}\}x)]y  && \text{by (\Bol)}\,.
\end{alignat*}
Thus using (\Bol) twice, we obtain
\begin{align*}
R((xy)^2) &= R([y(\{x[(y^{-1}x^{-1})(yx)]^{-1}\}x)]y) \\
&= R(y)R(\{x[(y^{-1}x^{-1})(yx)]^{-1}\}x)R(y) \\
&= R(y)R(x)R([(y^{-1}x^{-1})(yx)]^{-1})R(x)R(y)\,,
\end{align*}
as claimed.
\end{proof}

\section{Proof of \thmref{thm:main}}
\label{sec:results}

In this section, we give the proof of \thmref{thm:main}.

\medskip

1) Assume that for any $x,y\in L$, either $xy=yx$ or $x^{-1}(xy)=y$.
Applying this hypothesis to $x$ and $[y(xy)^{-1}]y$, we have that either
\begin{equation}
\label{eqn:hmm2}
x \{[y(xy)^{-1}]y\} = \{[y(xy)^{-1}]y\} x
\qquad\text{or}\qquad
x^{-1}(x\{[y(xy)^{-1}]y\}) = [y(xy)^{-1}]y\,.
\end{equation}
Now
\begin{equation}
\label{eqn:hmm3}
x \{[y(xy)^{-1}]y\} = [(xy)(xy)^{-1}]y = y
\end{equation}
using (\Bol), and so the left side of the left option of \eqref{eqn:hmm2} reduces to $y$.
For the right side, we use \eqref{eqn:soln} (with the roles of $x$ and $y$
reversed) to get $\{[y(xy)^{-1}]y\} x = \{[x^{-1}(yx)]x^{-1}\} x = x^{-1}(yx)$,
using (\RIP). Thus the left equation in \eqref{eqn:hmm2} reduces to
$x^{-1}(yx) = y$. Applying \eqref{eqn:hmm3}, the right option of \eqref{eqn:hmm2}
reduces to $x^{-1}y = [y(xy)^{-1}]y$. Cancelling $y$'s and using (\RIP), this
becomes $x^{-1}(xy) = y$. Summarizing, \eqref{eqn:hmm2} simplifies to
\[
x^{-1}(yx) = y
\qquad\text{or}\qquad
x^{-1}(xy) = y\,.
\]
Combining this with our hypothesis, we have that for all $x,y\in L$, either
\[
(\quad xy = yx\quad\text{and}\quad x^{-1}(yx) = y\quad )
\qquad\text{or}\qquad
x^{-1}(xy) = y\,.
\]
We see that in either case, the left inverse property $x^{-1}(xy) = y$ holds,
and so $L$ is Moufang by \propref{thm1}.

\medskip

2) Assume that for any $x,y\in L$,
either $xy=yx$ or $(xy)x = x(yx)$.  Applying this hypothesis to $x^{-1}$
and $(xy)x$ we have that either
\[
x^{-1}(xy\cdot x)=(xy\cdot x)x^{-1}
\qquad\text{or}\qquad
[x^{-1}(xy\cdot x)]x^{-1}=x^{-1}[(xy\cdot x)x^{-1}]\,.
\]
To the left option, we apply (\Bol) to the left side and (\RIP) to the
right side to get $yx = xy$. To the right option, we apply (\Bol) followed
by (\RIP) to the left side and (\RIP) to the right side to get
$y=x^{-1}(xy)$. Thus either $xy = yx$ or $x^{-1}(xy) = y$.  By part~(1),
$L$ is Moufang.

\medskip

3) Assume that for any $x,y\in L$, either
$x^{-1}(xy)=y$ or $(xy)^{-1}=y^{-1}x^{-1}$.
Apply the hypothesis to $x^{-1}$ and $xy$ to get
\[
x[x^{-1}(xy)] = xy
\qquad\text{or}\qquad
[x^{-1}(xy)]^{-1} = (xy)^{-1} x\,.
\]
The left option reduces to $x^{-1}(xy) = y$, and so for any $x,y\in L$,
we have either
\[
x^{-1}(xy) = y
\qquad\text{or}\qquad
(\quad (xy)^{-1} = y^{-1}x^{-1}
\quad\text{and}\quad
[x^{-1}(xy)]^{-1} = (xy)^{-1} x \quad )\,.
\]
The right option implies $[x^{-1}(xy)]^{-1} = (y^{-1}x^{-1})x = y^{-1}$ by (\RIP),
and so $x^{-1}(xy) = y$. Thus in either case, the left inverse property
$x^{-1}(xy) = y$ holds identically, and so $L$ is Moufang by \propref{thm1}.

\medskip

4) Assume that for any $x,y\in L$, either
$(xy)x=x(yx)$ or $(yx)^{-1}=x^{-1}y^{-1}$.
Applying this hypothesis to $x^{-1}$ and $(xy)x$, we get
\[
[x^{-1}((xy)x)]x^{-1} = x^{-1}[((xy)x)x^{-1}]
\qquad\text{or}\qquad
[((xy)x)x^{-1}]^{-1} = x ((xy)x)^{-1}\,.
\]
The left side of the left option reduces to $y$,
using (\Bol) and (\RIP). The right side reduces
to $x^{-1}(xy)$ using (\RIP). Thus the left option
simplifies to $x^{-1}(xy) = y$. The left side of
the right option reduces to $(xy)^{-1}$ using (\RIP).
Using (\SAIP) and then (\Bol), the right side becomes
$x((x^{-1}y^{-1})x^{-1}) = y^{-1}x^{-1}$. Thus the right
option simplifies to $(xy)^{-1} = y^{-1}x^{-1}$.
Summarizing, we have for each $x,y\in L$, either
$x^{-1}(xy) = y$ or $(xy)^{-1} = y^{-1}x^{-1}$.
By part (3), $L$ is Moufang.

\medskip

5) Assume that for any $x,y\in L$, either $x(x^{-1}y) = y$ or $y^{-1}(yx) = x$,
where we have already used \lemref{lip-equiv} in the left equation.
Apply this to $x^{-1}$ and $yx$ to get
\[
x^{-1}(x\cdot yx)=yx
\qquad\text{or}\qquad
(yx)^{-1}(yx\cdot x^{-1}) = x^{-1}\,.
\]
We rewrite the left option as
$x^{-1}(x\cdot yx)=yx=(x^{-1}x\cdot y)x = x^{-1}(xy\cdot x)$, using (\Bol)
and so cancelling, we have $x(yx) = (xy)x$. The right option simplifies to
$(yx)^{-1}y = x^{-1}$ using (\RIP) and, using it once more,
we have $(yx)^{-1} = x^{-1}y^{-1}$.
Summarizing, we have that for all $x,y\in L$,
$x(yx) = (xy)x$ or $(yx)^{-1} = x^{-1}y^{-1}$. By Part (4), $L$ is Moufang.

\medskip

Next we show the equivalence of conditions (6), (7) and (8).

\smallskip

First, assume (6), that is, for any $x,y\in L$, either
$(xy)x = x(yx)$ or $(xy)^{-1} = y^{-1}x^{-1}$. Apply this
hypothesis to $x^{-1}$ and $(xy)x$ to get
\[
[x^{-1}((xy)x)]x^{-1} = x^{-1}[((xy)x)x^{-1}]
\qquad\text{or}\qquad
[x^{-1}((xy)x)]^{-1}=((xy)x)^{-1} x\,.
\]
The left side of the left option simplifies to $y$
using (\Bol) and (\RIP). The right side reduces to
$x^{-1}(xy)$ using (\RIP). Thus the left option
reduces to $x^{-1}(xy) = y$. The left side of the
right option reduces to $(yx)^{-1}$ by (\Bol), while the
right side reduces to $x^{-1}y^{-1}$ using (\SAIP)
and then (\RIP). Thus the right option reduces to
$(yx)^{-1} = x^{-1}y^{-1}$. Summarizing, (6) implies
that for all $x,y\in L$, either $x^{-1}(xy) = y$ or
$(yx)^{-1} = x^{-1}y^{-1}$. This is precisely (7).

Next, assume (7), that is, for any $x,y\in L$,
$x(x^{-1}y)=y$ or $(yx)^{-1} = x^{-1}y^{-1}$
where we have already used \lemref{lip-equiv} in the left equation.
Apply this to $x^{-1}$ and $yx$ to get
\[
x^{-1}[x(yx)] = yx
\qquad\text{or}\qquad
((yx)x^{-1})^{-1} = x(yx)^{-1}\,.
\]
Since $yx = [(x^{-1}x)y]x = x^{-1}[(xy)x]$ using (\Bol),
the left option can written as $x^{-1}[x(yx)] = x^{-1}[(xy)x]$. Canceling,
we obtain $x(yx) = (xy)x$. The left side of the right option reduces
$y^{-1}$ by (\RIP), and so the right option becomes $y^{-1}=x(yx)^{-1}$,
which is equivalent to $y^{-1}(yx) = x$ by (\RIP). Summarizing, (7) implies
that for all $x,y\in L$, $(xy)x = x(yx)$ or $y^{-1}(yx) = x$.
This is condition (8).

Finally, assume (8), that is, for any $x,y\in L$, either
$(xy)x = x(yx)$ or $y^{-1}(yx) = x$. Apply this $x$ and $xy$
to get
\[
[x(xy)] x = x[(xy)x]
\qquad\text{or}\qquad
(xy)^{-1}[(xy)x] = x\,.
\]
For the left option, we have
$[x(xy)]x=x[(xy)x] = (x^2 y) x$ by (\Bol), and so canceling gives us
$x(xy) = x^2 y$. By \corref{cor:1}, $(xy)x = x(yx)$. Next, using
(\RIP), we may rewrite the
right equation as $(xy)^{-1} = x[(xy)x]^{-1} =
x[(x^{-1}y^{-1})x^{-1}] = y^{-1}x^{-1}$, where we have used (\SAIP) and
(\Bol) in the second and third equalities, respectively.
Summarizing, (8) implies that for all $x,y\in L$,
$(xy)x = x(yx)$ or $(xy)^{-1} = y^{-1}x^{-1}$, which is precisely
condition (6).

\bigskip

Now we are ready to show that each of conditions (6), (7) and (8)
implies that $L$ is Moufang, and by the preceding discussion, we
may assume that all three conditions hold.
We show that $L$ is Moufang by showing that $L$ satisfies the
flexible law identically. To this end, we assume that there exist
elements $a,b\in L$ such that $(ab)a \neq a(ba)$ or equivalently
(by \corref{cor:1}) $a^2 b \neq a(ab)$. By (6), we have
\begin{equation}
\label{eqn:ab-aaip}
(ab)^{-1} = b^{-1} a^{-1}\,.
\end{equation}
By (8) and \lemref{lip-equiv}, we have the following:
\begin{equation}
\label{eqn:ba-lip2}
b (b^{-1}a) = a\,.
\end{equation}

In (7), we take $x = b$ and $y = a^{-1}$ to obtain
\[
b(b^{-1}a^{-1}) = a^{-1}
\qquad\text{or}\qquad
(a^{-1}b)^{-1} = b^{-1}a\,,
\]
where we have already applied \lemref{lip-equiv} to the left equation.
That equation becomes $a^{-1} = b(ab)^{-1}$ using \eqref{eqn:ab-aaip},
or $a^{-1}(ab) = b$ using (\RIP). The right equation is equivalent to
$1 = (a^{-1}b)(b^{-1}a)$, which implies $b = [(a^{-1}b)(b^{-1}a)]b
= a^{-1}\{[b(b^{-1}a)]b\} = a^{-1}(ab)$ using (\Bol) followed by
\eqref{eqn:ba-lip2}. Thus in either case, we may use \lemref{lip-equiv}
to conclude:
\begin{equation}
\label{eqn:ab-lip2}
a (a^{-1}b) = b\,.
\end{equation}

Applying (\RIP) to \eqref{eqn:ab-lip2}, we have
$b(a^{-1}b)^{-1} = a = b(b^{-1}a)$ by \eqref{eqn:ba-lip2}.
Canceling, we get
\begin{align}
(a^{-1} b)^{-1} &= b^{-1} a \label{eqn:aib-aaip}
\intertext{and similarly,}
(b^{-1} a)^{-1} &= a^{-1} b\,. \label{eqn:bia-aaip}
\end{align}

Next, in (6) we take $x = a$ and $y = ab$ to get
\[
[a(ab)]a = a[(ab)a]
\qquad\text{or}\qquad
[a(ab)]^{-1} = (ab)^{-1} a^{-1}\,.
\]
The right side of the left equation is $(a^2 b) a$ by (\Bol),
and canceling $a$'s, we then have $a(ab) = a^2 b$. This is a
contradiction, and so the right equation must hold. Its right
side is $(b^{-1}a^{-1}) a^{-1} = b^{-1} a^{-2}$ using
\eqref{eqn:ab-aaip} and (\RAP). Thus we have
\begin{equation}
\label{eqn:ab-weird}
[a(ab)]^{-1} = b^{-1} a^{-2}\,.
\end{equation}

Next, in (7), take $x = b$ and $y = a^2$ to get
\[
b^{-1}(ba^2) = a^2
\qquad\text{or}\qquad
(a^2 b)^{-1} = b^{-1} a^{-2}\,.
\]
The right option implies $(a^2 b)^{-1} = [a(ab)]^{-1}$
by \eqref{eqn:ab-weird}, or $a^2 b = a(ab)$, a contradiction.
Hence the left option must hold, and so using \lemref{lip-equiv},
we have
\begin{equation}
\label{eqn:weird2}
b(b^{-1} a^2) = a^2\,.
\end{equation}

Now in (8), take $x = a$ and $y = a^{-1} b$ to get
\[
a[a(a^{-1}b)] = a^2 (a^{-1}b)
\qquad\text{or}\qquad
(a^{-1}b)^{-1}[(a^{-1}b)a] = a\,,
\]
where we have already applied \corref{cor:1} to the
left equation. That equation then reduces to
$ab = a^2 (a^{-1}b)$ by \eqref{eqn:ab-lip2}, and thus
$(ab)a = [a^2 (a^{-1}b)]a = a\{[a(a^{-1}b)]a\}
= a(ba)$ by (\Bol) and \eqref{eqn:ab-lip2}. This is
a contradiction, and so the right equation must hold.
Using \eqref{eqn:aib-aaip}, we have
$a = (b^{-1}a)[(a^{-1}b)a]$, and so
$b^{-1} a = a[(a^{-1}b)a]^{-1}$. Thus by (\RAP),
$b^{-1} a^2 = (b^{-1} a)a = \{a[(a^{-1}b)a]^{-1}\}a$.
Now we apply \eqref{eqn:soln} to this last expression
to get $b^{-1} a^2 = \{ (a^{-1} b)^{-1} [a(a^{-1}b)]\} (a^{-1} b)^{-1}$.
By \eqref{eqn:ab-lip2} and \eqref{eqn:aib-aaip}, this becomes
$b^{-1} a^2 = \{ (b^{-1}a)b\} (b^{-1}a)$. Now multiply on the left by $b$
and use \eqref{eqn:weird2} to get
$a^2 = b(b^{-1}a^2) = b [\{ (b^{-1}a)b\} (b^{-1}a)]$. By (\Bol) and
\eqref{eqn:ba-lip2}, this is equivalent to
$a^2 = \{[b(b^{-1}a)]b\} (b^{-1}a) = (ab)(b^{-1}a)$. Multiplying on
the right by $b$, we have
$a^2 b = [(ab)(b^{-1}a)]b = a\{ [b(b^{-1}a)] b\} = a(ab)$, where
we used (\Bol) and \eqref{eqn:ba-lip2} in the second and third
equalities, respectively. This last contradiction proves the desired
result.

\medskip

9) Assume that for all $x,y\in L$, either $(xx)y = x(xy)$ or
$x^{-1}(xy) = y$, where we have already used \corref{cor:1}. It will
be useful to write this as follows: for all $x,y\in L$, either
\begin{equation}
\label{eqn:case9}
x R(x,y) = x \qquad\text{or}\qquad x^{-1} R(x,y) = x^{-1}\,.
\end{equation}
We will verify the left inverse property, and so by way of contradiction,
assume there exist $a,b\in L$ such that
$a^{-1}(ab) \neq b$, that is, $a^{-1} R(a,b) \neq a^{-1}$.
By hypothesis, $a^2 b = a(ab)$, that is, $a R(a,b) = a$.

Take $x = b(ab)^{-1}$ and $y = b$ in \eqref{eqn:case9}, and use
$R(b(ab)^{-1},b) = R((ab)^{-1},b)$ (by \eqref{eqn:rim-ls}) to get
\[
[b(ab)^{-1}] R((ab)^{-1},b) = b(ab)^{-1}
\quad\text{or}\quad
[b(ab)^{-1}]^{-1} R((ab)^{-1},b) = [b(ab)^{-1}]^{-1}\,.
\]
Now $[b(ab)^{-1}]^{-1} = [(ab)(ab)^{-1}][b(ab)^{-1}]^{-1}
= a R(b) R((ab)^{-1}) R(b(ab)^{-1})^{-1} = a R(b,(ab)^{-1})$. Thus
the right option becomes
$[b(ab)^{-1}]^{-1} = a R(b,(ab)^{-1}) R((ab)^{-1},b) = a$ by
\eqref{eqn:rim-inv}. Hence this option implies $b(ab)^{-1} = a^{-1}$,
that is, $b = a^{-1}(ab)$ by (\RIP). This is a contradiction, and so the
left option must hold, which we record as:
\begin{equation}
\label{eqn:case9b}
[b(ab)^{-1}] R((ab)^{-1},b) = b(ab)^{-1}\,.
\end{equation}

We now compute
\begin{align*}
R(a,(ab)^{-1} b) &= R((ab)b^{-1},(ab)^{-1} b) && \text{by (\RIP)} \\
&= R(b^{-1},ab) R((ab)^{-1},b) && \text{by \lemref{lem:decomp}} \\
&= R((ab)b^{-1},ab) R((ab)^{-1},b) && \text{by \eqref{eqn:rim-ls}} \\
&= R(a,ab) R((ab)^{-1},b) && \text{by (\RIP)} \\
&= R(a,b) R((ab)^{-1},b) && \text{by \eqref{eqn:rim-rs}}\,.
\end{align*}
So taking $x = a$ and $y = (ab)^{-1} b$ in \eqref{eqn:case9}, we then have
\[
a R(a,b) R((ab)^{-1},b) = a
\qquad\text{or}\qquad
a^{-1} R(a,b) R((ab)^{-1},b) = a^{-1}\,.
\]
Since $a^{-1} R(a,b) = b(ab)^{-1}$, the right option simplifies to
$b(ab)^{-1} = a^{-1}$, using \eqref{eqn:case9b}. Thus $b = a^{-1}(ab)$
by (\RIP), a contradiction. Hence the left option holds. It simplifies
to $a R((ab)^{-1},b) = a$, since $a R(a,b) = a$ by the discussion following
\eqref{eqn:case9}. By \eqref{eqn:rim-inv},
$a = a R(b,(ab)^{-1}) = [b(ab)^{-1}]^{-1}$, and so as before,
$b(ab)^{-1} = a^{-1}$, that is, $b = a^{-1}(ab)$ by (\RIP). This final
contradiction completes the proof of this case.

\medskip

10) Assume that for all $x,y\in L$,
\begin{equation}
\label{eqn:ff}
(xy)x = x(yx) \qquad\text{or}\qquad (yx)y = y(xy)\,.
\end{equation}
We will verify the conditions of case (8). To this end, assume there
exist $a,b\in L$ such that $(ab)a \neq a(ba)$ and $b^{-1}(ba)\neq a$.
By hypothesis,
\begin{equation}
\label{eqn:b-flex}
(ba)b = b(ab)\,.
\end{equation}

We compute
\begin{align*}
b[a(ab)] &= b\{[(ab)b^{-1}](ab)\} && \text{by (\RIP)} \\
&= \{[b(ab)]b^{-1}\}(ab) && \text{by (\Bol)} \\
&= \{[(ba)b]b^{-1}\}(ab) && \text{by \eqref{eqn:b-flex}} \\
&= (ba)(ab) && \text{by (\RIP)}\,,
\end{align*}
and we record the result of this calculation as
\begin{equation}
\label{eqn:114}
b[a(ab)] = (ba)(ab)\,.
\end{equation}

In \eqref{eqn:ff}, replace $x$ with $x^{-1}$ and then replace $y$ with $yx$.
By (\RIP), the left option becomes $[x^{-1}(yx)]x^{-1} = x^{-1}y$, or, again by (\RIP),
$x^{-1}(yx) = (x^{-1}y)x$. By \lemref{lip-equiv},
this is equivalent to $x^{-1}(xy) = y$. Switching $x$ and $y$ in
\corref{cor:1}, the right option of \eqref{eqn:ff} may be written
in the form $y(yx) = y^2 x$. Again replacing $x$ by $x^{-1}$ and
$y$ by $yx$, this becomes $(yx)^2 x^{-1} = (yx)[(yx)x^{-1}] = (yx)y$
by (\RIP). Again by (\RIP), this is
equivalent to $yx = (yx)R(yx)R(x)^{-1}R(y)^{-1} =
(yx)R(y,x)^{-1} = (yx)R(x,y)$, using \eqref{eqn:rim-inv}. Thus
$(yx)(xy) = (yx)R(x)R(y) = (yx^2)y$ by (\RAP). Summarizing, we
have that for each $x,y\in L$, either
\begin{equation}
\label{eqn:ff-star1}
x^{-1}(xy) = y \qquad\text{or}\qquad (yx^2)y = (yx)(xy)\,.
\end{equation}
Since $b^{-1}(ba)\neq a$, we also have
\begin{equation}
\label{eqn:ff-br}
(ab^2)a = (ab)(ba)\,.
\end{equation}

Next, in \eqref{eqn:ff-star1}, replace $x$ with $x^{-1}$ and then
replace $y$ with $(xy)x$. The left option becomes
$x\{x^{-1}[(xy)x]\} = (xy)x$ or $x(yx) = (xy)x$ using (\Bol).
The left side of the right option becomes
\begin{align*}
\{[(xy)x]x^{-2}\}[(xy)x] &= [(\{[(xy)x]x^{-1}\})x^{-1}][(xy)x] && \text{by (\RAP)} \\
&= [(xy)x^{-1}][(xy)x] && \text{by (\RIP)} \\
&= (\{[(xy)x^{-1}]x\}y)x  && \text{by (\Bol)} \\
&= [(xy)y]x && \text{by (\RIP)} \\
&= (xy^2)x && \text{by (\RAP)}\,.
\end{align*}
The right side of the right option becomes
$\{[(xy)x]x^{-1}\}\{x^{-1}[(xy)x]\} = (xy)(yx)$, using (\RIP) and (\Bol).
Summarizing, we have that for each $x,y\in L$, either
\begin{equation}
\label{eqn:ff-star2}
x(yx) = (xy)x \qquad\text{or}\qquad (xy^2)x = (xy)(yx)\,.
\end{equation}

In \eqref{eqn:ff}, take $x = ab$ and $y = a$. The right option
becomes $[a(ab)]a = a[(ab)a] = (a^2 b)a$ by (\Bol) or $a(ab) = a^2 b$. By
\corref{cor:1}, this contradicts our assumption, and so the left option
holds, that is, $[(ab)a](ab) = (ab)[a(ab)]$. By \corref{cor:1}, this
is equivalent to $(ab)^2 a = (ab)[(ab)a)] = \{[(ab)a]b\}a
= \{a[(ba)b]\}a$, using (\Bol) twice. Thus,
\begin{equation}
\label{eqn:142}
(ab)^2 = a[(ba)b]\,.
\end{equation}

Next, in \eqref{eqn:ff-star1}, take $x = b^{-1}$ and $y = a$. The left
option becomes $b(b^{-1}a) = a$, which is a contradiction by
\lemref{lip-equiv}. Thus the right option holds, that is,
$(ab^{-2})a = (ab^{-1})(b^{-1}a)$. Using (\RAP) and rearranging, this
is equivalent to $(ab^{-1})R(b^{-1})R(a)R(b^{-1}a)^{-1} = ab^{-1}$,
that is, $(ab^{-1})R(b^{-1},a) = ab^{-1}$. By \eqref{eqn:rim-inv},
we have $ab^{-1} = (ab^{-1})R(a,b^{-1}) = \{[(ab^{-1})a]b^{-1}\} R(ab^{-1})^{-1}
= \{a[(b^{-1}a)b^{-1}]\} R(ab^{-1})^{-1}$ by (\Bol). Hence,
\begin{equation}
\label{eqn:143}
(ab^{-1})^2 = a[(b^{-1}a)b^{-1}]\,.
\end{equation}

Next, in \eqref{eqn:ff-star2}, take $x = b^{-1}$ and $y = ab$. Note that
$yx = a$ by (\RIP). The left option becomes $b^{-1}a = [b^{-1}(ab)]b^{-1}$
or $(b^{-1}a)b = b^{-1}(ab)$ by (\RIP). By \lemref{lip-equiv}, this is
equivalent to $b^{-1}(ba) = a$, a contradiction. Therefore the right option
holds, that is,
\begin{align*}
[b^{-1}(ab)]a &= [b^{-1}(ab)^2]b^{-1} && \\
&= b^{-1} R((ab)^2) R(b^{-1}) && \\
&= b^{-1} R(b) R(a) R([(b^{-1}a^{-1})(ba)]^{-1})R(a)R(b)R(b^{-1})
&& \text{by \eqref{eqn:turnsquare}} \\
&= \{ a[(b^{-1}a^{-1})(ba)]^{-1}\}a\,.
\end{align*}
Canceling $a$'s, we have
\begin{align*}
b^{-1}(ab) &= a[(b^{-1}a^{-1})(ba)]^{-1} && \\
&= a [(\{[(ba)^{-1}(ba)][(b^{-1}a^{-1})(ba)]^{-1}\}(ba))(ba)^{-1}] && \text{by (\RIP)} \\
&= a(\{(ba)^{-1}[\{(ba)[(b^{-1}a^{-1})(ba)]^{-1}\}(ba)]\}(ba)^{-1}) && \text{by (\Bol)} \\
&= a(\{(ba)^{-1}[(ab^2)a]\}(ba)^{-1}) && \text{by \eqref{eqn:preaip}} \\
&= a(\{(ba)^{-1}[(ab)(ba)]\}(ba)^{-1}) && \text{by \eqref{eqn:ff-br}} \\
&= a(\{(ab)[(ba)(ab)]^{-1}\}(ab)) && \text{by \eqref{eqn:soln}} \\
&= \{[a(ab)][(ba)(ab)]^{-1}\}(ab) && \text{by (\Bol)}\,.
\end{align*}
Canceling $ab$, we have
$b^{-1} = [a(ab)][(ba)(ab)]^{-1}$ or $a(ab) = b^{-1}[(ba)(ab)]$ by (\RIP). Finally,
applying \eqref{eqn:114}, we obtain
\begin{equation}
\label{eqn:201}
b^{-1}\{b[a(ab)]\} = a(ab)\,.
\end{equation}

Now suppose $(ba^2)b = b[a(ab)]$. Then by \eqref{eqn:201} and (\Bol), we have
$a(ab) = b^{-1}\{b[a(ab)]\} = b^{-1}[(ba^2)b] = a^2 b$.
This contradicts our assumption, and so we conclude
\begin{equation}
\label{eqn:contra}
(ba^2)b \neq b[a(ab)]\,.
\end{equation}

Next, in \eqref{eqn:ff}, with the left option in the form $x^2y = x(xy)$ (\corref{cor:1}),
take $x = ba$ and $y = a^{-1}$. Note that $xy = b$ by (\RIP). The left option becomes
$(ba)^2 a^{-1} = (ba)b = b(ab)$ by \eqref{eqn:b-flex}. Thus
\begin{align*}
b &= [(ba)^2 a^{-1}] (ab)^{-1} && \text{by (\RIP)} \\
&= (ba) R(ba)R(a)^{-1}R(b)^{-1} R(b) R((ab)^{-1}) && \\
&= (ba) R(b,a)^{-1} R(b) R((ab)^{-1}) && \\
&= (ba) R(a,b) R(b) R((ab)^{-1}) && \text{by \eqref{eqn:rim-inv}} \\
&= \{[(ba)a]b\} R((ab)^{-1})R(b)R((ab)^{-1}) && \text{by (\RIP)} \\
&= \{[(ba)a]b\} R([(ab)^{-1}b](ab)^{-1}) && \text{by (\Bol)} \\
&= \{[(ba)a]b\} R([(ab)b^{-1}](ab))^{-1} && \text{by (\SAIP) and (\RIP)} \\
&= \{[(ba)a]b\} R(a(ab))^{-1} && \text{by (\RIP)} \\
&= [(ba^2)b][a(ab)]^{-1} && \text{by (\RAP) and (\RIP)}\,.
\end{align*}
By (\RIP), $b[a(ab)] = (ba^2)b$, but this contradicts \eqref{eqn:contra}.
Therefore the right option holds of \eqref{eqn:ff}, with $x = ba$ and $y = a^{-1}$,
holds, that is, $a^{-1}b = [a^{-1}(ba)]a^{-1}$,
or $(a^{-1}b)a = a^{-1}(ba)$ by (\RIP). By \lemref{lip-equiv}, this is
equivalent to
\begin{equation}
\label{eqn:206}
a^{-1}(ab) = b\,.
\end{equation}

Next, in \eqref{eqn:ff}, with the left option in the form $x^2y = x(xy)$ (\corref{cor:1}),
take $x = a^{-1}(ba)$ and $y = (ba)^{-1}$. Note that $xy = a^{-1}$ by (\RIP).
The left option is $[a^{-1}(ba)]a^{-1} = [a^{-1}(ba)]^2 (ba)^{-1}$, and so
by (\RIP),
$a^{-1}(ba) = [a^{-1}(ba)] R(a^{-1}(ba))R(ba)^{-1}R(a) =
[a^{-1}(ba)] R(a^{-1},ba)^{-1}$. By \eqref{eqn:rim-inv},
$a^{-1}(ba) = [a^{-1}(ba)] R(ba,a^{-1})
= [a^{-1} (ba)^2] R(a)^{-1} R(b)^{-1}$ by (\RAP) and (\RIP). Hence
\begin{align*}
a^{-1}(ba)^2 &= \{[a^{-1}(ba)]b\}a && \\
&= [(\{[a^{-1}(ba)]a^{-1}\}a)b]a && \text{by (\RIP)} \\
&= [(\{[b(ab)^{-1}]b\}a)b]a && \text{by \eqref{eqn:soln}} \\
&= \{[b(ab)^{-1}][(ba)b]\}a && \text{by (\Bol)} \\
&= (\{[b(ab)](ab)^{-2}\}[(ba)b])a && \text{by (\RAP) and (\RIP)} \\
&= (\{[b(ab)]\{a[(ba)b]\}^{-1}\}[(ba)b])a && \text{by \eqref{eqn:142}} \\
&= [(a^{-1}\{[(ba)b]a\})a^{-1}]a && \text{by \eqref{eqn:soln}} \\
&= a^{-1}\{[(ba)b]a\} && \text{by (\RIP)}\,.
\end{align*}
Canceling $a^{-1}$'s, we have $(ba)^2 = [(ba)b]a$. Thus
$ba = (ba)R(b)R(a)R(ba)^{-1} = (ba)R(b,a)$, and so
$(ba)R(a,b) = ba$ by \eqref{eqn:rim-inv}. Thus $[(ba)a]b = (ba)(ab)$.
Apply (\RAP) to the left side and \eqref{eqn:114} to the right side, we have
$(ba^2)b = b[a(ab)]$. This contradicts \eqref{eqn:contra}. We conclude that
the right option of \eqref{eqn:ff}, with $x = a^{-1}(ba)$ and $y = (ba)^{-1}$,
must hold, and this is
\begin{align*}
(ba)^{-1} a^{-1} &= \{(ba)^{-1}[a^{-1}(ba)]\}(ba)^{-1} && \\
&= (\{(ba)[a^{-1}(ba)]^{-1}\}(ba))^{-1} && \text{by (\SAIP)} \\
&= (\{a[(ba)a^{-1}]\}a)^{-1} && \text{by \eqref{eqn:soln}} \\
&= (a^{-1}b^{-1})a^{-1} && \text{by (\RIP) and (\SAIP)}\,.
\end{align*}
Canceling $a^{-1}$'s, we obtain
\begin{equation}
\label{eqn:207}
(ba)^{-1} = a^{-1} b^{-1}\,.
\end{equation}

Next, in \eqref{eqn:ff-star1}, take $x = a$ and $y = b^{-1}$. The left option
is $a^{-1}(ab^{-1}) = b^{-1}$ and so by \lemref{lip-equiv} and \eqref{eqn:207},
$b^{-1} = a(a^{-1}b^{-1}) = a(ba)^{-1}$. By (\RIP), $b^{-1}(ba) = a$, a
contradiction. Thus the right option holds, that is,
$(b^{-1}a^2)b^{-1} = (b^{-1}a)(ab^{-1})$. By (\RAP) and rearranging,
this is equivalent to
$b^{-1}a = (b^{-1}a)R(a)R(b^{-1})R(ab^{-1})^{-1} = (b^{-1}a)R(a,b^{-1})$.
By \eqref{eqn:rim-inv}, $b^{-1}a = (b^{-1}a)R(b^{-1},a)$ and so
$(b^{-1}a)^2 = [(b^{-1}a)b^{-1}]a$. Applying (\Bol), we obtain
\begin{equation}
\label{eqn:212}
(b^{-1}a)^2 = b^{-1}[(ab^{-1})a]\,.
\end{equation}

Finally, in \eqref{eqn:ff-star1}, take $x = b^{-1}a$ and $y = [(ab^{-1})a]^{-1}
= (a^{-1}b)a^{-1}$ (using (\SAIP)). Note that
$xy = (b^{-1}a)[(a^{-1}b)a^{-1}] = \{[(b^{-1}a)a^{-1}]b\}a^{-1}
= (b^{-1}b)a^{-1} = a^{-1}$, using (\Bol) and (\RIP). The left option is
$(b^{-1}a)^{-1} a^{-1} = (a^{-1}b)a^{-1}$, or $(b^{-1}a)^{-1} = a^{-1}b$
after canceling. By \eqref{eqn:206} (and \lemref{lip-equiv}),
$b = a(a^{-1}b) = a(b^{-1}a)^{-1}$, and so $b(b^{-1}a) = a$ by (\RIP).
By \lemref{lip-equiv}, this is equivalent to $b^{-1}(ba) = a$, a contradiction.
Thus the right option must hold, that is,
\begin{equation}
\label{eqn:lastright}
\{[(a^{-1}b)a^{-1}](b^{-1}a)^2\}[(a^{-1}b)a^{-1}] = \{[(a^{-1}b)a^{-1}](b^{-1}a)\}a^{-1}\,.
\end{equation}
By (\SAIP), the left side of \eqref{eqn:lastright} is
\begin{align*}
\{[(ab^{-1})a]^{-1}(b^{-1}a)^2\}[(ab^{-1})a]^{-1} &=
\{[(ab^{-1})a]^{-1}(b^{-1}[(ab^{-1})a])\}[(ab^{-1})a]^{-1} && \text{by \eqref{eqn:212}} \\
&= (b^{-1}\{[(ab^{-1})a]b^{-1}\}^{-1})b^{-1} && \text{by \eqref{eqn:soln}} \\
&= (b^{-1}\{a[(b^{-1}a)b^{-1}]\}^{-1})b^{-1} && \text{by (\Bol)} \\
&= [b^{-1}(ab^{-1})^{-2}]b^{-1} && \text{by \eqref{eqn:143}}\,.
\end{align*}
By (\Bol), the right side of \eqref{eqn:lastright} is
\begin{align*}
(a^{-1}b)\{[a^{-1}(b^{-1}a)]a^{-1}\} &= [((a^{-1}b)\{[a^{-1}(b^{-1}a)]a^{-1}\})b]b^{-1}
&& \text{by (\RIP)} \\
&= \{a^{-1}[(b\{[a^{-1}(b^{-1}a)]a^{-1}\})b]\}b^{-1} && \text{by (\Bol)} \\
&= \{a^{-1} [\{[(ba^{-1})(b^{-1}a)]a^{-1}\} b] \} b^{-1} && \text{by (\Bol)} \\
&= \{a^{-1} [(ba^{-1})R(b^{-1}a,a^{-1})]\} b^{-1} && \text{by (\RIP)}\,.
\end{align*}
Now
\[
R(b^{-1}a,a^{-1}) = R(b^{-1}a,(b^{-1}a)a^{-1})
= R(b^{-1}a,b^{-1}) = R(a,b^{-1})\,,
\]
using \eqref{eqn:rim-rs}, (\RIP) and \eqref{eqn:rim-ls}. Continuing, the right side of
\eqref{eqn:lastright} is now
\begin{align*}
\{a^{-1} [(ba^{-1})R(a,b^{-1})]\} b^{-1}
&= \{ a^{-1} [(ba^{-1})R(a)R(b^{-1})R(ab^{-1})^{-1}]\}b^{-1} && \\
&= \{ a^{-1} (ab^{-1})^{-1}\} b^{-1} && \text{by (\RIP)}\,.
\end{align*}
Putting this together, we have reduced \eqref{eqn:lastright} to
\[
[b^{-1}(ab^{-1})^{-2}]b^{-1} = \{ a^{-1} (ab^{-1})^{-1}\} b^{-1}\,.
\]
Canceling $b^{-1}$'s, we have
$a^{-1}(ab^{-1})^{-1} = b^{-1}(ab^{-1})^{-2} = [b^{-1}(ab^{-1})^{-1}](ab^{-1})^{-1}$
by (\RAP). Canceling once more, we obtain
$a^{-1} = b^{-1}(a b^{-1})^{-1}$, and so $b^{-1} = a^{-1}(ab^{-1})$. By
\lemref{lip-equiv} and \eqref{eqn:207}, $b^{-1} = a(a^{-1}b^{-1}) = a(ba)^{-1}$.
By (\RIP) once more, we get $b^{-1}(ba) = a$, our final contradiction.

We have shown that with \eqref{eqn:ff} as hypothesis,
the assumption that there exist elements $a,b\in L$ satisfying
$(ab)a \neq a(ba)$ and $b^{-1}(ba) \neq a$
is untenable. It follows that for all $x,y\in L$, either $(xy)x = x(yx)$ or
$y^{-1}(yx) = x$. This is case (8), and so it follows
that $L$ is a Moufang loop.

\bigskip

This completes the proof of \thmref{thm:main}.

\section{Right Automorphic Bol Loops}
\label{sec:auto}

In the Introduction, we admitted that we can say nothing about the
imposition of one particular pair of equations for elements $x,y$
in a Bol loop $L$, namely:
\[
\text{Either}\qquad (xy)^{-1} = y^{-1}x^{-1}
\qquad\text{or}\qquad (yx)^{-1} = x^{-1}y^{-1}\,.  \tag{AA}
\]
In this section and the next, we will show that (AA) is
sufficient for a Bol loop to be Moufang in two special classes of
Bol loops.

In this section, we consider the class of \emph{right automorphic} Bol loops,
that is, Bol loops $L$ in which each
right inner mapping $R(x,y)$ is an automorphism:
\[
(uv) R(x,y) = [uR(x,y)] [vR(x,y)]  \tag{A$_r$}
\]
for all $x,y,u,v\in L$. This class includes other interesting classes of
Bol loops as special cases.
A \emph{Bruck loop} is a Bol loop with the automorphic inverse property
$(xy)^{-1} = x^{-1}y^{-1}$. A \emph{Burn loop} is a right conjugacy closed Bol loop,
that is, a Bol loop in which the set of right translations is closed under conjugation:
for all $x,y$, there exists $z$ such that $R(x)R(y)R(x)^{-1} = R(z)$.
Both Bruck loops and Burn loops are right automorphic. For Bruck loops, see, \emph{e.g.},
\cite{Kiechle:02} in the dual setting of left Bol loops and with
different terminology. That Burn loops are right automorphic follows from the more
general fact that right conjugacy closed loops have this property; see \cite{CC}

The right automorphic property can also be expressed in terms of right translations:
\[
R(z) R(x,y) = R(x,y) R(zR(x,y))\,.  \tag{A$_r$}
\]
We will also use the fact that automorphisms commute with inversion, that is, (A$_r$)
implies
\begin{equation}
\label{eqn:inv-comm}
[x R(y,z)]^{-1} = x^{-1} R(y,z)
\end{equation}
for all $x,y,z$.

\begin{lem}
\label{lem:rim-ar}
For all $x,y$ in a right automorphic Bol loop,
\begin{equation}
\label{eqn:rim-ar}
R(x^{-1},y^{-1}) = R(x,y)\,.
\end{equation}
\end{lem}

\begin{proof}
First observe that, by (\RIP), $R(x^{-1}y^{-1})^{-1} =
R((x^{-1}y^{-1})^{-1})$, so that, with
$z = (x^{-1}y^{-1})^{-1}$ and the roles of $x$ and $y$ interchanged,
(A$_r$) becomes
\begin{align*}
R((x^{-1}y^{-1})^{-1}) R(y,x)
&= R(y,x) R((x^{-1}y^{-1})^{-1} R(y,x)) && \\
&= R(y,x) R([(x^{-1}y^{-1}) R(y,x)]^{-1}) && \text{by \eqref{eqn:inv-comm}} \\
&= R(y,x) R((yx)^{-1})^{-1} && \text{by (\RIP)} \\
&= R(y,x) R(yx) && \text{by (\RIP)} \\
&= R(y) R(x)\,.
\end{align*}
Therefore,
\begin{align*}
R(x^{-1},y^{-1}) R(y,x) &=
R(x)^{-1} R(y)^{-1} R((x^{-1}y^{-1})^{-1}) R(y,x) && \text{by (\RIP)} \\
&= R(x)^{-1} R(y)^{-1} R(y) R(x) && \\
&= \mathrm{id}_L\,.
\end{align*}
Applying \eqref{eqn:rim-inv}, we obtain the desired result.
\end{proof}

\begin{thm}
\label{thm:bol-ar}
Let $L$ be a right automorphic Bol loop satisfying (AA). Then $L$ is a Moufang loop.
\end{thm}

\begin{proof}
Assume that $x,y\in L$ satisfy
\begin{equation}
\label{eqn:xy-aaip}
(xy)^{-1} = y^{-1} x^{-1}\,.
\end{equation}
We will show that this implies $x^{-1}(xy) = y$. First, observe that
$[x(xy)^{-1}] R(xy,x^{-1}) = ((xy)x^{-1})^{-1}$. Now
\begin{align*}
R(xy,x^{-1}) &= R((xy)^{-1},x) && \text{by \eqref{eqn:rim-ar}} \\
&= R(y^{-1}x^{-1},x) && \text{by \eqref{eqn:xy-aaip}} \\
&= R(y^{-1}x^{-1},(y^{-1}x^{-1})x) && \text{by \eqref{eqn:rim-rs}} \\
&= R(y^{-1}x^{-1},y^{-1}) && \text{by (\RIP)} \\
&= R(x^{-1},y^{-1}) && \text{by \eqref{eqn:rim-ls}} \\
&= R(x,y) && \text{by \eqref{eqn:rim-ar}}\,.
\end{align*}
Thus $((xy)x^{-1})^{-1} = [x(xy)^{-1}] R(x,y) = [x R(x,y)] [((xy)^{-1}) R(x,y)]
= [x R(x,y)][ (y^{-1}x^{-1}) R(x,y)] = [x R(x,y)] (xy)^{-1}$ using (A$_r$)
in the second equality, \eqref{eqn:xy-aaip} in the third, and then simplifying,
using (\RIP). We record this as
\begin{equation}
\label{eqn:ar-trick}
((xy)x^{-1})^{-1} = [x R(x,y)] (xy)^{-1}\,.
\end{equation}

Now apply (AA) to $xy$ and $x^{-1}$ to get that either
\[
((xy)x^{-1})^{-1} = x(xy)^{-1} \qquad\text{or}\qquad
(x^{-1}(xy))^{-1} = (xy)^{-1} x \,.
\]
The right option reduces to $(x^{-1}(xy))^{-1} = (xy)^{-1} x = (y^{-1}x^{-1})x = y^{-1}$
using \eqref{eqn:xy-aaip} and (\RIP). Inverting both sides, we get $x^{-1}(xy) = y$.
For the left option, we use \eqref{eqn:ar-trick} to get
$[x R(x,y)](xy)^{-1} = x(xy)^{-1}$, or just $x R(x,y) = x$ after canceling.
Inverting both sides of this and using \eqref{eqn:inv-comm},
we have $x^{-1} R(x,y) = x^{-1}$, which, after rearranging, is equivalent
to $x^{-1}(xy) = y$. Summarizing, both options imply $x^{-1}(xy) = y$.

Returning once more to (AA), we see that for $x,y\in L$, the left option,
which is $(xy)^{-1} = y^{-1}x^{-1}$, implies $x^{-1}(xy) = y$.
Thus $L$ satisfies the condition of case (7) of
\thmref{thm:main}. Therefore $L$ is a Moufang loop.
\end{proof}

\section{SRAR Loops}
\label{sec:srar}

In this section we consider another important class of Bol loops for which
condition (AA) implies that the Moufang identity holds, namely those with
(strongly) right alternative loop rings.

For any commutative, associative ring $R$ with unity and any loop $L$,
one forms the \emph{loop ring} $RL$ just as if $L$ were a group.
Loop rings that are not associative but which satisfy interesting identities
are hard to find, but they do exist.
Here is an example of relevance to this paper.

\begin{prop}[\cite{EGG:95a}]
\label{SRARtheorem}
If $L$ is a loop and $R$ is a commutative ring with unity and with
$\ch R=2$, then the loop ring $RL$ satisfies the right Bol identity
if and only if $L$ is a (right) Bol loop and, for all $x,y,z,w\in L$,
at least one of the following three conditions holds:
\begin{align*}
\tag*{$D(x,y,z,w)$:} [(xy)z]w &= x[(yz)w] \quad\text{and}\quad [(xw)z]y=x[(wz)y]   \\
\tag*{$E(x,y,z,w)$:} [(xy)z]w &= x[(wz)y] \quad\text{and}\quad [(xw)z]y=x[(yz)w]  \\
\tag*{$F(x,y,z,w)$:} [(xy)z]w &= [(xw)z]y \quad\text{and}\quad x[(yz)w]=x[(wz)y]\,.
\end{align*}
\end{prop}

For simplicity, we will often refer to the above conditions by saying that
$D$ or $E$ or $F$ holds, and we will follow a similar convention for other
conditions.

A ring (with unity) satisfying the Bol identity necessarily satisfies the right
alternative law, but the converse is not true. Thus
rings satisfying the Bol identity are called \emph{strongly right alternative}.
A (necessarily Bol) loop $L$ is called \emph{SRAR}, for \emph{strongly right
alternative ring}, if, for some commutative, associative ring $R$ with unity,
$RL$ is strongly right alternative but not (left) alternative. We refer
the interested reader to \cite{EGG:95a,EGG:96a,EGG:05d,EGG:06c,EGG:07a,EGG:08a}
for the relevant literature on SRAR loops (whose existence is well established).

It is worth recording that under the assumption that the underlying loop $L$
is Bol, the conditions of \propref{SRARtheorem}
can be considerably simplified. Of course, in the right equation in
condition $F(x,y,z,w)$, the $x$'s can be canceled. In addition, the left equation
in condition $D(x,y,z,w)$ can be written in the form
$x R(y)R(z)R(w) R((yz)w)^{-1} = x$. By \lemref{lem.5}
with $\phi = R(y)R(z)R(w)$, this is
equivalent to $x = x R(w)R(z)R(y) R((wz)y)^{-1}$, but this, in turn,
is equivalent to the right equation in condition $D(x,y,z,w)$.
So only one equation is required in condition $D(x,y,z,w)$.

Next, if we distribute the disjunctions over the conjunctions, we have
that a Bol loop $L$ has a strongly right alternative loop ring if and
only if it satisfies each of the following four conditions for all $x,y,z,w\in L$:
\begin{alignat}{8}
[(xy)z]w &= x[(yz)w] &\quad\text{or}\quad [(xy)z]w &= x[(wz)y] &&\quad\text{or}\quad &[(xy)z]w &= [(xw)z]y\,, \label{eqn:def1}\\
[(xy)z]w &= x[(yz)w] &\quad\text{or}\quad [(xy)z]w &= x[(wz)y] &&\quad\text{or}\quad &(yz)w &= (wz)y\,, \label{eqn:def2}\\
[(xw)z]y &= x[(wz)y] &\quad\text{or}\quad [(xw)z]y &= x[(yz)w] &&\quad\text{or}\quad &[(xy)z]w &= [(xw)z]y\,, \label{eqn:def3}\\
[(xw)z]y &= x[(wz)y] &\quad\text{or}\quad [(xw)z]y &= x[(yz)w] &&\quad\text{or}\quad &(yz)w &= (wz)y\,. \label{eqn:def4}
\end{alignat}
Here we have used the equivalence of the equations in $D(x,y,z,w)$ in the first options of these
conditions. But now it is evident that \eqref{eqn:def1} and \eqref{eqn:def3} are equivalent
by the symmetry in $y$ and $w$, and similarly, \eqref{eqn:def2} and \eqref{eqn:def4} are equivalent.
Discarding \eqref{eqn:def3} and \eqref{eqn:def4}, and rearranging what remains, we have shown the
following.

\begin{thm}
\label{shortSRARtheorem}
A loop $L$ has a strongly right alternative loop ring if and only if
$L$ is a Bol loop and for all $x,y,z,w\in L$, at least one of the following
conditions holds:
\begin{align*}
\tag*{$D_0(x,y,z,w)$} [(xy)z]w &= x[(yz)w]\,, \\
\tag*{$E_0(x,y,z,w)$} [(xy)z]w &= x[(wz)y]\,, \\
\tag*{$F_0(x,y,z,w)$} [(xy)z]w = [(xw)z]y \quad &\text{and}\quad (yz)w = (wz)y\,.
\end{align*}
\end{thm}

If we set $w=1$ (or $y=1$ or $z=1$) in the conditions $D$, $E$ and $F$,
we get conditions that characterize the loop ring being right alternative .

\begin{prop}[\cite{EGG:88a}]
\label{prp:rar}
A loop $L$ has a right alternative loop ring if and only if
$L$ is right alternative, and for all $x,y,z\in L$, at least one of the
following conditions holds:
\begin{align*}
\tag*{$D'(x,y,z)$:} (xy)z &= x(yz) \quad\text{and}\quad (xz)y=x(zy)\,, \\
\tag*{$E'(x,y,z)$:} (xy)z &= x(zy) \quad\text{and}\quad (xz)y=x(yz)\,, \\
\tag*{$F'(x,y,z)$:} (xy)z &= (xz)y \quad\text{and}\quad x(yz)=x(zy)\,.
\end{align*}
\end{prop}

In the Bol case, these conditions simplify even further. In the same way
that we derived \thmref{shortSRARtheorem} from Proposition \ref{SRARtheorem}, so
from Proposition \ref{prp:rar} we can derive the following.

\begin{cor}\label{cor:shortRAR}
A Bol loop $L$ has a right alternative loop ring if and only if,
for all $x,y,z\in L$,  at least one of the following conditions holds:
\begin{align*}
\tag*{$D'_0(x,y,z)$} (xy)z &= x(yz)\,, \\
\tag*{$E'_0(x,y,z)$} (xy)z &= x(zy)\,, \\
\tag*{$F'_0(x,y,z)$} (xy)z = (xz)y \quad &\text{and}\quad yz = zy\,.
\end{align*}
\end{cor}

Now consider ``$D'_0(x^{-1},x,y)$ or $E'_0(x^{-1},x,y)$ or $F'_0(x^{-1},x,y)$''.
This gives us for each $x,y$ in a Bol loop $L$ with a right alternative loop ring:
\begin{equation}
\label{eqn:rarweak-pre}
x^{-1}(xy) = y \quad\text{or}\quad x^{-1}(yx) = y  \quad\text{or}\quad
( xy = yx \quad\text{and}\quad x^{-1}y = yx^{-1} )\,.
\end{equation}
Now obviously we may weaken this still further as: for each $x,y\in L$,
\begin{equation}
\label{eqn:rarweak}
x^{-1}(xy) = y \quad\text{or}\quad x^{-1}(yx) = y  \quad\text{or}\quad xy = yx\,.
\end{equation}
Conversely, suppose that \eqref{eqn:rarweak} holds where $L$ is a Bol loop.
Replace $x$ with $x^{-1}$ to obtain
\begin{equation}
\label{eqn:rar-weak-pre2}
x(x^{-1}y) = y \quad\text{or}\quad x(yx^{-1}) = y  \quad\text{or}\quad x^{-1}y = yx^{-1}\,.
\end{equation}
By \lemref{lip-equiv}, the first option is equivalent to $x^{-1}(xy) = y$.
Multiply both sides of the second option on the right by $x$ and on the left
by $x^{-1}$. The left side then becomes
$x^{-1}\{[x(yx^{-1})]x\} = (yx^{-1})x = y$ by (\Bol) and (\RIP), and so
the second option is now $y = x^{-1}(yx)$. Thus \eqref{eqn:rar-weak-pre2}
is now
\[
x^{-1}(xy) = y \quad\text{or}\quad x^{-1}(yx) = y  \quad\text{or}\quad x^{-1}y = yx^{-1}\,.
\]
Putting this together with \eqref{eqn:rarweak}, we have that \eqref{eqn:rarweak-pre}
holds. So we may now refer to \eqref{eqn:rarweak} instead of \eqref{eqn:rarweak-pre}
without loss of generality.

There do exist (nonMoufang) Bol loops with right alternative loop rings which are not SRAR,
and there exist (nonMoufang) Bol loops satisfying \eqref{eqn:rarweak} which do not have
right alternative loop rings. Using the LOOPS package \cite{LOOPS} for GAP \cite{GAP},
we found that of the $2033$ nonMoufang Bol loops of order $16$, $1873$ are SRAR loops,
$5$ have a right alternative loop ring that is not strongly right alternative, and
$2$ satisfy \eqref{eqn:rarweak} but do not have a right alternative loop ring.

We do not know if \eqref{eqn:rarweak} itself has an interpretation in terms of loop rings,
that is, \emph{if} $L$ \emph{is a Bol loop such that for all} $x,y\in L$,
\emph{\eqref{eqn:rarweak} holds, then what, if anything, can one say about the loop ring}
$RL$\emph{?}

The connection between SRAR loops (and the other classes of Bol loops mentioned
in this section) and the theme of this paper appears in the next theorem.

\begin{thm}\label{thm2} Let $L$ be a Bol loop such that for all $x,y\in L$,
\eqref{eqn:rarweak} holds. Suppose also that for each $x,y\in L$,
\[
\text{either}\qquad (xy)^{-1}=y^{-1}x^{-1}
\qquad\text{or}\qquad
(yx)^{-1}=x^{-1}y^{-1}\,.  \tag{AA}
\]
Then $L$ is Moufang.  In particular, no SRAR loop (which, by definition is not
Moufang) can satisfy (AA).
\end{thm}

\begin{proof}
We will show that the hypotheses imply condition (3) of \thmref{thm:main}.
To this end, assume there exist $a,b\in L$ such that $a^{-1}(ab)\neq b$
and $(ab)^{-1} \neq b^{-1} a^{-1}$. By \eqref{eqn:rarweak},
\begin{equation}
\label{eqn:weak-ab}
a^{-1}(ba) = b \qquad\text{or}\qquad ab = ba\,.
\end{equation}
By (AA),
\begin{equation}
\label{eqn:ba-aaip}
(ba)^{-1} = a^{-1} b^{-1}\,.
\end{equation}

Of the alternatives in \eqref{eqn:weak-ab}, let us assume first that
$a^{-1}(ba) = b$. Then $b^{-1} = [a^{-1}(ba)]^{-1}$ and so
$(ba)b^{-1} = (ba)[a^{-1}(ba)]^{-1} = [(aa^{-1})(ba)][a^{-1}(ba)]^{-1}
= a R(a^{-1},ba)$. Now $R(a^{-1},ba) = R((ba)a^{-1},ba) = R(b,ba) = R(b,a)$,
using \eqref{eqn:rim-ls}, (\RIP), and \eqref{eqn:rim-rs}. Thus
$(ba)b^{-1} = a R(b,a) = [(ab)a] R(ba)^{-1}$, and so
$(ab)a = [(ba)b^{-1}](ba)$.
Inverting both sides of this and using (\SAIP) on the right side, we have
$[(ab) a]^{-1} = [(ba)^{-1}b](ba)^{-1} =
[(a^{-1} b^{-1})b](ba)^{-1} = a^{-1}(ba)^{-1}$,
where we have used \eqref{eqn:ba-aaip} and (\RIP) in the second and third
equalities, respectively. Thus
\begin{equation}
\label{eqn:aba-mess}
(ab)a = [a^{-1}(ba)^{-1}]^{-1}\,.
\end{equation}

Now by (AA) with $x = a^{-1}$ and $y = (ba)^{-1}$,
\[
[a^{-1}(ba)^{-1}]^{-1} = (ba)a \qquad\text{or}\qquad
[(ba)^{-1} a^{-1}]^{-1} = a(ba)\,.
\]
With the left option, \eqref{eqn:aba-mess} gives $(ab)a = (ba)a$ and so $ab = ba$.
For the right option, we use \eqref{eqn:ba-aaip} and (\SAIP) to get
$a(ba) = [(a^{-1}b^{-1})a^{-1}]^{-1} = (ab)a$.

Summarizing what we have so far, we have shown that \emph{if} $a^{-1}(ba) = b$,
\emph{then} $ab = ba$ \emph{or} $a(ba) = (ab)a$. Thus in \eqref{eqn:weak-ab}, we
may replace $a^{-1}(ba) = b$ with $a(ba) = (ab)a$, that is, we have
\begin{equation}
\label{eqn:newopts1}
a(ba) = (ab)a \qquad\text{or}\qquad ab = ba\,.
\end{equation}

We continue with the assumption that $a^{-1}(ba) = b$, but first we make
some observations using the fact that, in a loop, the solution of an
equation is unique. That is, if $xy = z$ and $xw = z$, then $y = w$.

Note that, by (\Bol), $a^{-1}\{[b(a^{-1}b)^{-1}]b\}
= [(a^{-1}b)(a^{-1}b)^{-1}]b = b$, so that $[b(a^{-1}b)^{-1}]b$ is a solution
of $a^{-1}x = b$. But, by our assumption, $a^{-1}(ba) = b$, so that $ba$
is also a solution. Therefore, $[b(a^{-1}b)^{-1}]b = ba$, or, by (\RIP),
$b(a^{-1}b)^{-1} = (ba)b^{-1}$.

Also by (\Bol),
$[(ba)b^{-1}][((ba)\{[(ba)b^{-1}](ba)\}^{-1})(ba)]
= (\{[(ba)b^{-1}](ba)\} \{[(ba)b^{-1}](ba)\}^{-1})(ba)
= ba$, so that $y = ((ba)\{[(ba)b^{-1}](ba)\}^{-1})(ba)$
is a solution of $[(ba)b^{-1}] y = ba$. But, as we saw above,
$[b(a^{-1}b)^{-1}]b = ba$, so that $b$ is also a solution.
Hence, by the uniqueness of solutions,
$((ba)\{[(ba)b^{-1}](ba)\}^{-1})(ba) = b = a^{-1}(ba)$, by
our assumption. Canceling the $ba$'s, we conclude that
$(ba)\{[(ba)b^{-1}](ba)\}^{-1} = a^{-1}$, and so
$ba = a^{-1}\{[(ba)b^{-1}](ba)\}$ by (\RIP).
Thus $z = [(ba)b^{-1}](ba)$ is a solution of $a^{-1}z = ba$,
as is $(ab)a$, the latter by (\Bol). Once again using the
uniqueness of solutions, $[(ba)b^{-1}](ba) = (ab)a$.

Summarizing, we have shown that \emph{if} $a^{-1}(ba) = b$, \emph{then}
$[b(a^{-1}b)^{-1}](ba) = (ab)a$. Now, if the left option of
\eqref{eqn:newopts1} holds, then $a(ba) = (ab)a =
[b(a^{-1}b)^{-1}](ba)$. Canceling the $ba$'s, we have $a = (ba)b^{-1}$,
so that $ab = ba$ by (\RIP).

Thus, even if the left options hold in both \eqref{eqn:weak-ab} and
\eqref{eqn:newopts1}, we still get $ab = ba$, so that, in any case,
\begin{equation}
\label{eqn:ab-comm}
ab = ba\,.
\end{equation}

Now since (\Bol) and \eqref{eqn:ab-comm} imply
\begin{equation}
\label{eqn:home}
a = b[(a(ba)^{-1})a] = b[(a(ab)^{-1})a]\,,
\end{equation}
we have
\begin{align*}
a^{-1} b^{-1} &= (ba)^{-1} && \text{by \eqref{eqn:ba-aaip}} \\
&= (ab)^{-1} && \text{by \eqref{eqn:ab-comm}}\\
&= (\{b[(a(ab)^{-1})a]\}b)^{-1} && \text{by \eqref{eqn:home}} \\
&= (b^{-1} [(a^{-1}(ab))a^{-1}]) b^{-1} &&\text{by (\SAIP), twice} \\
&= \{[(b^{-1} a^{-1})(ab)]a^{-1}\} b^{-1} &&\text{by (\Bol)}\,.
\end{align*}
Canceling twice, we have $1 = (b^{-1} a^{-1})(ab)$, that is, $(ab)^{-1} = b^{-1} a^{-1}$.
This is our final contradiction.

We have shown that the assumptions that both $a^{-1}(ab) \neq b$ and $(ab)^{-1} \neq b^{-1} a^{-1}$
are untenable. It follows that for all $x,y\in L$, either $x^{-1}(xy) = y$ or
$(xy)^{-1} = y^{-1}x^{-1}$. This is case (3) of \thmref{thm:main}, and so it follows
that $L$ is a Moufang loop.
\end{proof}

Years ago, the first two authors showed that if $L$ is Moufang and
$RL$ satisfies the right alternative law, then $RL$ is actually an
alternative ring \cite{EGG:88a}.  So this corollary is immediate.

\begin{cor}
If $L$ is a Bol loop with a right alternative loop ring $RL$,
and if for any $x,y\in L$, (AA) holds, then $RL$ is an alternative ring.
\end{cor}

\section{More on the DEF Conditions}

The previous section has focussed renewed attention on conditions $D'$, $E'$ and $F'$,
at least one of which must hold in any SRAR loop.  In fact, all three
conditions are essential in the sense that if $L$ is a Bol
loop and one of $D',E',F'$ is discarded---that is, for all $x,y,z\in L$,
one of the remaining two conditions is satisfied---then $L$ is not SRAR: it is a
(possibly associative) Moufang loop.

\begin{thm}\label{thm:dpepfp}
Let $L$ be a Bol loop. If, for each $x,y,z\in L$,
\begin{enumerate}
\item $D'(x,y,z)$ or $E'(x,y,z)$ holds, then $L$ is a Moufang loop.
\item $D'(x,y,z)$ or $F'(x,y,z)$ holds, then $L$ is a group,
\item $E'(x,y,z)$ or $F'(x,y,z)$ holds, then $L$ is an abelian group.
\end{enumerate}
\end{thm}
\begin{proof}
1) Suppose for any $x,y,z\in L$, either  $D'(x,y,z)$ or $E'(x,y,z)$ holds.
Then, in particular,
$D'(x^{-1},xy,x)$ or $E'(x^{-1},xy,x)$ for each $x,y\in L$.
Thus, for all $x,y\in L$,
\begin{equation}
\label{main-3.5} [x^{-1}(xy)]x = x^{-1}[(xy)x] \qquad\text{or}\qquad
xy=(x^{-1}x)xy = x^{-1}[(xy)x]\,.
\end{equation}
By (\Bol), $x^{-1}[(xy)x] =[(x^{-1}x)y]x=
yx$. Using this in both parts of \eqref{main-3.5}, and then
canceling $x$'s in the left equation, we obtain
\begin{equation*}
x^{-1}(xy) = y \qquad\text{or}\qquad xy = yx
\end{equation*}
for all $x,y\in L$. So  $L$ is Moufang by Part~(1) of \thmref{thm:main}.
\medskip\par
2) We first verify that $L$ is a Moufang loop. If,
for each $x,y,z\in L$, $D'(x,y,z)$ or $F'(x,y,z)$ holds, then for
each $x,y\in L$, $D'(x,x^{-1},y)$ or $F'(x,x^{-1},y)$ holds.
Thus,
\begin{displaymath}
y = x(x^{-1}y) \qquad\text{or}\qquad y=(xy)x^{-1}.
\end{displaymath}
Applying \lemref{lip-equiv} to the first equation here and (\RIP)
to the second, we obtain
\begin{displaymath}
y = x^{-1}(xy) \qquad\text{or}\qquad xy = yx.
\end{displaymath}
This holds for all $x,y\in L$, so $L$ is Moufang by
Part~(1) of \thmref{thm:main}.

Now, we show that $L$ is associative.  Let $x,y,z\in L$.
It suffices to show that $D'(x,y,z)$ holds.  The contrary gives
us $F'(x,y,z)$ and, applying Moufang's Theorem, we also have $F'(z,x,y)$,
$F'(y,x,z)$ and $F'(z,yz,x)$.
From the second equation of $F'(z,x,y)$, we get
\begin{equation}
xy=yx;
\end{equation}
from the first equation of $F'(y,x,z)$, we get
\begin{equation}
(yx)z=(yz)x;
\end{equation}
and from the second equation of $F'(z,yz,x)$, we get
\begin{equation}
(yz)x=x(yz).
\end{equation}
Putting these together, we have
\begin{equation}
(xy)z=(yx)z=(yz)x=x(yz),
\end{equation}
the desired result.
\medskip\par
3) Taking $x=1$ in $E'(x,y,z)$ and $F'(x,y,z)$, we see immediately
that $L$
is commutative and hence Moufang.  Take $x,y,z\in L$. Either
$(xy)z=x(zy)=x(yz)$ by $E'(x,y,z)$ and commutativity, or
$z(xy)=(xy)z=(xz)y=(zx)y$ by commutativity and $F'(x,y,z)$.  In
either case, $x$, $y$ and $z$ associate. Thus $L$ is an abelian
group.
\end{proof}

\begin{cor}\label{notodd}
A Bol loop of odd order with a strongly right alternative
loop ring is a group.
\end{cor}

\begin{proof}
We first show that any Bol loop of odd order with a right alternative
loop ring must be Moufang. Either $D'(x,y,z)$ or $E'(x,y,z)$ or $F'(x,y,z)$
holds for any $x,y,z\in L$.

Suppose $F'(x,y,z)$ for some $x,y,z$; that is, $(xy)z = (xz)y$ and $yz=zy$.
Then $x R(y)R(z) R(yz)^{-1} = x R(z)R(y) R(zy)^{-1}$, that is,
$x R(y,z) = x R(z,y)$. By \eqref{eqn:rim-inv},
$x R(y,z)^2 = x$.

Now let $G(L)$ denote the group generated by the right translations of $L$,
that is, $G(L) = \langle R_x\ |\ x\in L\rangle$. By a result essentially due
to Glauberman (\cite[Theorem 14]{Gl}), given explicitly in
\cite[Theorem 4.14]{FKP}, $G(L)$ has odd order. Thus $R(x,y)$, which is an
element of $G(L)$, has odd order. It follows that $x R(y,z) = x$, that is,
$(xy)z = x(yz)$, so that $D'(x,y,z)$ holds.

Summarizing, we have that for any triple $x, y, z\in L$, either $D'(x,y,z)$
must hold or $E'(x,y,z)$ holds.
But then, by \thmref{thm:dpepfp}, $L$ is Moufang.
It follows that $RL$ is alternative \cite[Theorem 1.3]{EGG:88a},
and so $L$ is associative \cite[Corollary 2.5]{EGG:90a}.
\end{proof}

\providecommand{\MR}{\relax\ifhmode\unskip\space\fi MR }
\providecommand{\MRhref}[2]{%
  \href{http://www.ams.org/mathscinet-getitem?mr=#1}{#2}
}
\providecommand{\href}[2]{#2}

\end{document}